\documentclass[reqno,draft]{amsart}
\usepackage{amsfonts}%
\usepackage{amssymb}
\usepackage{mathrsfs}%
\usepackage{cite}
\usepackage{color}
\allowdisplaybreaks
\date{\today}
\usepackage{mathrsfs}
\usepackage{amssymb}
\usepackage{amsmath}
\usepackage{amsthm}
\usepackage{latexsym}
\usepackage{graphicx}

 \numberwithin{equation}{section}

\newcommand{\dv}{\mathrm{div}\,}

\newcommand{\R}{{\mathbb{R}^3}}

\newtheorem{Theorem}{Theorem}[section]
\newtheorem{Lemma}[Theorem]{Lemma}

\newtheorem{Remark}{Remark}[section]
\begin{document}

\title[ 3D full compressible Navier-Stokes system in a bounded domain]
 { Global strong solutions to the 3D full compressible Navier-Stokes system with vacuum in a bounded domain}

\author[J.-S. Fan]{Jishan Fan}
\address{ Department of Applied Mathematics,
 Nanjing Forestry University, Nanjing 210037, P.R.China}
\email{fanjishan@njfu.edu.cn}

\author[F.-C. Li]{Fucai Li}
\address{Department of Mathematics, Nanjing University, Nanjing
 210093, P.R. China}
 \email{fli@nju.edu.cn}
\begin{abstract}
In this short paper we establish the global well-posedness of strong solutions to the 3D full compressible Navier-Stokes system with vacuum in a bounded domain $\Omega\subset \mathbb{R}^3$ by the bootstrap argument  provided that the viscosity coefficients $\lambda$ and $\mu$ satisfy that $7\lambda>9\mu$
and the initial data $\rho_0$ and $u_0$ satisfy that $\|\rho_0\|_{L^\infty(\Omega)}$ and $\|\rho_0|u_0|^5\|_{L^1(\Omega)}$ are sufficient small.
 \end{abstract}


\maketitle

\section{Introduction}

In this short paper, we consider the following initial and boundary problem to the 3D full compressible Navier-Stokes system in a bounded domain $\Omega\subset \mathbb{R}^3$:
\begin{align}
&\partial_t\rho+\dv(\rho u)=0\ \ \mathrm{on}\ \ \Omega\times(0,\infty),\label{1.1}\\
&\partial_t(\rho u)+\dv(\rho u\otimes u)-\mu\Delta u-(\lambda+\mu)\nabla\dv u+\nabla p=0\ \ \mathrm{on}\ \ \Omega\times(0,\infty),\label{1.2}\\
&C_V\{\partial_t(\rho\theta)+\dv(\rho u\theta)\}-\Delta\theta+p\dv u\nonumber\\
& \qquad \qquad =\frac{\mu}{2}|\nabla u+\nabla u^t|^2+\lambda(\dv u)^2\ \ \mathrm{on}\ \ \Omega\times(0,\infty),\label{1.3}\\
&u=0,\ \frac{\partial\theta}{\partial n}=0\ \ \mathrm{on}\ \ \partial\Omega\times(0,\infty),\label{1.4}\\
&(\rho,\rho u,\rho\theta)(\cdot,0)=(\rho_0,\rho_0u_0,\rho_0\theta_0)\ \ \mathrm{in}\ \ \Omega.\label{1.5}
\end{align}
Here the unknowns $\rho,u,\theta$ denote the density, velocity and temperature of the fluid, respectively. The pressure $p:=R\rho\theta$ and the internal energy $e:=C_V\theta$ with positive constants $R$ and $C_V$. $\lambda$ and $\mu$ are two viscosity constants satisfying $\mu>0$ and $\lambda+\frac23\mu\geq0$.
$n$ is the unit outward normal vector to the smooth boundary $\partial\Omega$ of $\Omega$.

If the initial density $\rho_0$ has a positive lower bound, the global existence of small smooth solutions to the problem \eqref{1.1}--\eqref{1.5} was obtained in \cite{v,MN}  three decades ago.

If the initial data may contain  vaccum, Cho and Kim \cite{1} proved the local well-posedness of strong solutions to the problem \eqref{1.1}-\eqref{1.5} under some compatibility conditions:
\begin{align}
&-\mu\Delta u_0-(\lambda+\mu)\nabla\dv u_0+R\nabla(\rho_0\theta_0)=\sqrt{\rho_0}g_1,\label{1.6}\\
&\Delta\theta_0+\frac{\mu}{2}|\nabla u_0+\nabla u_0^t|^2+\lambda(\dv u_0)^2=\sqrt{\rho_0}g_2,\label{1.7}
\end{align}
with $g_1,g_2\in L^2(\Omega)$.

Recently, Huang and Li \cite{2} prove that the global well-posedness of strong solutions to the full compressible Navier-Stokes equations  in the  whole space $\R$
with smooth initial data which are of small energy but possibly large oscillations where the initial density is allowed to vanish, see also  \cite{WZ}.
However, the methods developed in \cite{2,WZ}  can not applied directly to  bounded domain case.

The aim of this paper is to prove that, although the initial density may contain vaccum,  the problem \eqref{1.1}-\eqref{1.5} still has a unique global strong solution for small
initial data.  Our results reads as
\begin{Theorem}\label{th1.1}
Let $0\leq\rho_0\in W^{1,6}(\Omega), u_0\in H_0^1(\Omega)\cap H^2(\Omega), 0\leq\theta_0\in H^2(\Omega)$ with $\frac{\partial\theta_0}{\partial n}=0$ on $\partial\Omega$ and \eqref{1.6}, \eqref{1.7} hold true. If
\begin{equation}
7\lambda>9\mu\ \ and\ \ \|\rho_0\|_{L^\infty}+\|\rho_0|u_0|^5\|_{L^1}\label{1.8}
\end{equation}
is  sufficient \ small, then the problem \eqref{1.1}-\eqref{1.5} has a unique global-in-time strong solution.
\end{Theorem}
\begin{Remark}\label{re1.1}
It is interesting to note that the initial temperature  $\theta_0$ need not be small in our results.
\end{Remark}
\begin{Remark}\label{re1.2}
It is possible to establish  a similar result for the full compressible magnetohydrodynamical  system.
\end{Remark}
\begin{Remark}
\label{re1.3} When $\Omega:=\R$ and consider the isentropic Navier-Stokes system, a similar result can be proved when $\|\rho_0\|_{L^p}$ is small for some large $p$ by the method developed here and a blow-up criterion
\begin{equation}
\rho\in L^p(\R\times(0,T))\label{1.9}
\end{equation}
proved in  \emph{\cite{3}}.
\end{Remark}
\begin{Remark}
\label{re1.4} A similar result holds true when the boundary condition $\frac{\partial\theta}{\partial n}=0$ on $\partial\Omega$ is replaced by $\theta=0$ on $\partial\Omega$.
\end{Remark}

To prove Theorem \ref{th1.1}, we will use the following abstract bootstrap argument or continuity argument (\!\cite{4}, Page 20).
\begin{Lemma} [\!\cite{4}]
\label{le1.1}
Let $T>0$. Assume that two statements $C(t)$ and $H(t)$ with $t\in[0,T]$ satisfy the following conditions:

(a) If $H(t)$ holds for some $t\in[0,T]$, then $C(t)$ holds for the same $t$;

(b) If $C(t)$ holds for some $t_0\in[0,T]$, then $H(t)$ holds for $t$ in a neighborhood of $t_0$;

(c) If $C(t)$ holds for $t_m\in[0,T]$ and $t_m\rightarrow t$, then $C(t)$ holds;

(d) $C(t)$ holds for at least one $t_1\in[0,T]$.

\noindent
Then $C(t)$ holds for all $t\in[0,T]$.
\end{Lemma}

We will also use the following regularity criterion (\!\cite{6}):
\begin{Lemma}
\label{le1.2}
Let $0\leq\rho_0\in W^{1,6}(\Omega), u_0\in H_0^1(\Omega)\cap H^2(\Omega), 0\leq\theta_0\in H^2(\Omega)$ with $\frac{\partial\theta_0}{\partial n}=0$ on $\partial\Omega$ and \eqref{1.6}, \eqref{1.7} hold true.
 If $\rho$ and $u$ satisfy
\begin{equation}
\rho\in L^\infty(\Omega\times(0,T))\ \ and\ \ u\in L^5(\Omega\times(0,T)),\label{1.10}
\end{equation}
then
\begin{align}
&\|\rho\|_{L^\infty(0,T;W^{1,6}(\Omega))}+\|u\|_{L^\infty(0,T;H^2(\Omega))}\nonumber\\
& \qquad \qquad\quad  +\|\theta\|_{L^\infty(0,T;H^2(\Omega))}+\|u\|_{L^2(0,T;W^{2,6}(\Omega))}\leq C_1.\label{1.11}
\end{align}
\end{Lemma}

The remainder of this paper is to the proof of Theorem \ref{th1.1}. Our proof is very short due to that it heavily depends on using Lemma \ref{le1.2}.

\section{Proof of Theorem \ref{th1.1}}

We will use the bootstrap argument and regularity criterion \eqref{1.10} to prove Theorem \ref{th1.1}.

Let $\delta>0$ be a fixed number, say
\begin{equation}
2\|\rho_0\|_{L^\infty}+2\|\rho_0|u_0|^5\|_{L^1}\leq \delta .\label{2.1}
\end{equation}
Denote by $H(t)$ the statement that, for $t\in [0,T]$,
\begin{equation}
\|\rho\|_{L^\infty(\Omega\times[0,t])}+\|u\|_{L^5(\Omega\times[0,t])}^5\leq\delta\label{2.2}
\end{equation}
and $C(t)$ the statement that
\begin{equation}
\|\rho\|_{L^\infty(\Omega\times[0,t])}+\|u\|_{L^5(\Omega\times[0,t])}^5\leq\frac\delta2.\label{2.3}
\end{equation}
The conditions (b)-(d) in Lemma \ref{le1.1} are clearly true and it remains to verify (a) under the condition \eqref{1.8}. Once this is verified, then the bootstrap argument would imply that $C(t)$, or \eqref{2.3} actually holds for any $t\in [0,T]$ and thus \eqref{1.11} holds true.

Now we assume that \eqref{2.2} holds true for some $t\in[0,T]$.
By Lemma \ref{le1.2}, we have
\begin{equation}
\|\rho\|_{L^\infty(0,t;W^{1,6})}+\|u\|_{L^\infty(0,t;H^2)}+\|\theta\|_{L^\infty(0,t;H^2)}+\|u\|_{L^2(0,t;W^{2,6})}\leq C_1.\label{2.4}
\end{equation}

Testing \eqref{1.1} by $\rho^{q-1}\ (q>2)$ and using \eqref{2.4}, we see that
$$  \frac{d}{dt}\|\rho\|_{L^q}\leq\left(1+\frac1q\right)\|\dv u\|_{L^\infty}\|\rho\|_{L^q}\leq C_0\|u\|_{W^{2,6}}\|\rho\|_{L^q},$$
which yields
$$\|\rho\|_{L^q}\leq\|\rho_0\|_{L^q}\exp\left(C_0\int_0^T\|u\|_{W^{2,6}} dt\right)\leq\|\rho_0\|_{L^q}\exp(C_0\sqrt T C_1).$$
Taking $q\rightarrow+\infty$ and letting $$\|\rho_0\|_{L^\infty}\ \ \mathrm{be\ \ small},$$ we arrive at
\begin{equation}
  \|\rho\|_{L^\infty(\Omega\times[0,t])}\leq\|\rho_0\|_{L^\infty}\exp(C_0\sqrt T C_1)\leq\frac\delta4.\label{2.5}
\end{equation}

When $7\mu>9\lambda$, we can adopt a technique of Hoff \cite{7} to bound the velocity in $L^5(\Omega\times[0,t])$ as follows.
Setting $q=5$ and testing \eqref{1.2} by $q|u|^{q-2}u$ and using \eqref{1.11}, we derive
\begin{align}
  &\frac{d}{dt}\int\rho|u|^q dx+\int\big\{q|u|^{q-2}[\mu|\nabla u|^2+(\lambda+\mu)(\dv u)^2+\mu(q-2)|\nabla|u||^2]\nonumber\\
  &\qquad \qquad \qquad \qquad  +q(\lambda+\mu)(\nabla|u|^{q-2})\cdot u\dv u\big\}dx\nonumber\\
  =&q\int p\dv(|u|^{q-2}u)dx=q\int R\rho\theta\dv(|u|^{q-2}u)dx\leq C_2\int\rho|u|^{q-2}|\nabla u|dx\nonumber\\
  \leq&\epsilon\int|u|^{q-2}|\nabla u|^2 dx+\frac{C_2^2}{\epsilon}\|\rho\|_{L^\infty}^2.\label{2.6}
\end{align}

Now, noticing that $|\nabla|u||\leq|\nabla u|$ and the condition $7\mu>9\lambda$, we have after a straight calculation that (\cite{8}):
\begin{align*}
&\mathrm{the\  second\  term\  on\  the\ left\ hand\   side \ of \eqref{2.6}}\\
\geq  & q\int|u|^{q-2}\left[\mu(q-1)-\frac{\lambda+\mu}{4}(q-2)^2\right]|\nabla u|^2 dx\\
  \geq & C_0\int|u|^{q-2}|\nabla u|^2 dx.
\end{align*}

Inserting the above inequality into \eqref{2.6} and taking $\epsilon=\frac{C_0}{2}$, we have $$\int\rho|u|^q dx+\frac{C_0}{2}\int_0^t\int|u|^{q-2}|\nabla u|^2dx ds\leq \int\rho_0|u_0|^q dx+\frac{C_2^2}{\epsilon}\|\rho\|_{L^\infty}^2T$$ which gives
\begin{align}
\int_0^t\int|u|^5 dx ds&\leq C_0\int_0^t\int|u|^{q-2}|\nabla u|^2 dx ds \nonumber\\
&\leq C_0\int\rho_0|u_0|^5 dx+\frac{C_2^2}{\epsilon}T\|\rho\|_{L^\infty}^2\nonumber\\
&\leq\frac\delta4.\label{2.7}
\end{align}

Summing up \eqref{2.5} and \eqref{2.7} gives \eqref{2.3}. Thus we have \eqref{1.10}.
Hence we can apply Lemma \ref{le1.2} to complete the proof of Theorem \ref{th1.1}.
\hfill$\square$



\medskip\medskip
{\bf Acknowledgements:}
 Fan is supported by NSFC (Grant No. 11171154). Li is supported partially by NSFC (Grant No. 11671193), Fundamental Research Funds for the Central
Universities (Grant No. 020314380014) and
   PAPD.


\end{document}